\documentclass{article}
\pdfoutput=1

\usepackage{pdfpages}
\usepackage{lmodern}
\usepackage[utf8]{inputenc}
\usepackage[T1]{fontenc}

\begin{document}
\title{Letters to Alan Weinstein about Courant algebroids}
\author{Pavol \v Severa}
\date{}
\maketitle
\section*{Introduction}
Between October 1998 and March 2000 I wrote several letters to Alan Weinstein about Courant algebroids (CAs). They still seem to be of some interest, so I decided to post them on the arXiv to make them more easily available.

Courant algebroids were introduced by Liu, Weinstein, and Xu \cite{lwx}; their principal motivation was extending the notion of Manin triples from Lie bialgebras to Lie bialgebroids. Similar structures were appearing in Poisson-Lie T-duality \cite{ks}, which was my original motivation for studying CAs.

Here is a short summary of the letters and some context.
\begin{description}
\item[Letter 1] contains the definition  and classification of exact CAs. It also gives a motivation for exact CAs coming from 2-dimensional variational problems, suggests a link between exact CAs and gerbes (motivated by the appearance of gerbes in 2-dim QFTs \cite{bry,gaw}), and contains a now-standard definition of CAs in terms of a non-skew-symmetric bracket. Some of these ideas were later put to \cite{SW}.
\item[Letter 2] is somewhat less relevant (perhaps besides linking Dirac structures with D-branes) - I tried there to explain  Poisson-Lie T-duality as it was my motivation for studying CAs. The proper link between CAs and Poisson-Lie T-duality is given  in Letter 4.
\item[Letter 3] contains local classification of transitive CAs and thoughts about global classification.
\item[Letter 4] contains a global classification of transitive CAs in terms of transitive Lie algebroids with vanishing 1st Pontryagin class, introduces reduction of CAs and formulates Poisson-Lie T-duality in terms of this reduction. It also introduces  examples of exact CAs over groups and homogeneous spaces. (Some of these ideas were later rediscovered, notably reduction of CAs in \cite{bur}, classification of transitive CAs in \cite{xu}, a link between CAs and T-duality \cite{cg}.)
\item[Letter 5] tries to make the link between exact CAs and gerbes more precise.
\item[Letter 6] gives a construction of a natural generating Dirac operator for any CA. It was written after I learned from Alan Weinstein that his student Dmitry Roytenberg, following ideas of Yvette Kosmann-Schwarzbach and of Arkady Vaintrob, found a link between (some) CAs and dg symplectic manifolds \cite{Rth}, that this link was improved by Weinstein himself who gave a construction of a dg symplectic manifold for any CA, and that there was a work in progress by Anton Alekseev and Ping Xu on generating Dirac operators for Courant algebroids \cite{ax}, which seemed a lot like a deformation quantization of those dg symplectic manifolds. This letter explains a natural version of such a deformation quantization.
\item[Letter 7] makes the link between  CAs and dg symplectic manifolds more precise - namely it states that CAs are equivalent to non-negatively graded  manifolds with a symplectic form of degree $2$ and with a degree-3 function $\theta$ satisfying $\{\theta,\theta\}=0$. As a part of this correspondence it gives a simple construction of a dg symplectic manifold, starting with a vector bundle $A\to M$ with a symmetric pairing, as a submanifold of $T^*[2]A[1]$. (Details were later published by Roytenberg \cite{Rst}.) In the first part of the letter the link between exact CAs and gebres is discussed from the dg point of view.
\item[Letter 8] suggests how to integrate CAs to symplectic 2-groupoids using ideas from Sullivan's rational homotopy theory \cite{Sul} and from the AKSZ construction \cite{AKSZ}. As a warm-up it gives a construction of a groupoid integrating a Lie algebroid by mimicking the construction of the fundamental groupoid of a topological space, which was apparently new at the time of writing. Last but not least, it introduces the silly terminology ``NQ-manifold''. (This letter was later expanded to \cite{some}.)
\end{description}



\section*{Supplementary material}

**Letter 1**

## How Courant algebroids appear in 2-dim variational problems (or maybe in string theory)

Suppose that we look for extremal surfaces in a manifold $M$. The thing we want to be extremal (stationary, more precisely) is the integral of a 2-form $\alpha$ living on $M$ (this is for simplicity, but after all it can always be achieved by passing to a jet space). According to Noether, if $v$ is a vector field preserving $\alpha$ (I shall use $v(\alpha)$ for the Lie derivative, i.e. $v(\alpha) = 0$) then $v \lrcorner \alpha$ is closed on extremal surfaces. Here it is just a consequence of Cartan's $v(\alpha) = d(v \lrcorner \alpha) + v \lrcorner d\alpha$. But conservation laws also come from symmetries up to total divergences. If $v$ is a vector field and $\theta$ a 1-form on $M$ and $v(\alpha) + d\theta = 0$ then $v \lrcorner \alpha + \theta$ is conserved (notice that we need $\theta$ and not just $d\theta$). We may say that the pair $(v, \theta)$ is a symmetry, but it would be nice to see it geometrically. And here the Courant algebroid structure on $A = TM \oplus T^*M$ appears. $(v, \theta)$ (from now on written as $v + \theta$) is a section of $A$. We also encode $\alpha$ in $A$ as a subbundle $C \subset A$; at $x \in M$ we simply take $C_x$ to be the graph of $\alpha_x : TM \to T^*M$. Now we'd like to see $v + \theta$ produce some flow on $A$ that would preserve $C$ iff $v(\alpha) + d\theta = 0$. This is pretty obvious: $v$ should give the natural flow on $A = TM \oplus T^*M$ while $\theta$ should move each fiber $A_x$ in itself with the velocity $-u \lrcorner d\theta$ at a point $u + \zeta \in A_x$. We may write down the corresponding bracket $[v_1 + \theta_1, v_2 + \theta_2]_A$ (the rate of change of $v_2 + \theta_2$ under the flow generated by $v_1 + \theta_1$):

$$[v_1 + \theta_1, v_2 + \theta_2]_A = [v_1, v_2] + v_1(\theta_2) - v_2 \lrcorner d\theta_1. \tag{1}$$

So the algebroid is here. Next we'll see that in fact only the algebroid structure on $A$ (including the projection $a : A \to TM$, the inner product $\langle v_1 + \theta_1, v_2 + \theta_2 \rangle = v_1 \lrcorner \theta_2 + v_2 \lrcorner \theta_1$ and the generation of flows by sections, i.e. the bracket $[,]_A$) is natural, while the explicit splitting $A = TM \oplus T^*M$ of the exact sequence $0 \to T^*M \to A \to TM \to 0$ is not (the second arrow in the sequence is $a^*$). The point is that $\alpha$ is not completely natural in variational problems: we can always add a closed $\beta$ (i.e. only $d\alpha$ matters; we can also start with a closed but not necessarily exact 3-form instead of $d\alpha$ – cf. the classification of exact Courant algebroids in terms of the 3rd deRham cohomology of $M$ below). What happens: adding a $\beta$ we produce an automorphism $v + \theta \mapsto v + \theta + v \lrcorner \beta$ of $A$. It is really an automorphism (this is immediately checked) of $A$ as a Courant algebroid, but it changes the splitting. The correct setting for 2-dim. var. problems is in fact this: an exact C. a. $A$ (one for which $0 \to T^*M \to A \to TM \to 0$ is exact – instead of $TM \oplus T^*M$) and an isotropic splitting $C$ of $A$ (instead of $\alpha$).

Here I enumerate the properties of C. a.'s, as I use the non-skewsymmetric bracket, that is natural here; it also has the virtue that the properties have clear meaning (except for one) and thus they are immediately checked. Most of them say that the flows generated by sections are automorphisms.

1. $[,]_A$ is really given by a flow, that is projected to a flow of $M$ (a vector field) simply by applying $a$ on the generating section of $A$:

$$[s_1, f s_2]_A = f[s_1, s_2] + a(s_1)(f) s_2$$

1

2. the flows preserve $a$:
$$a([s_1, s_2]_A) = [a(s_1), a(s_2)]$$

3. the flows preserve the inner product $\langle, \rangle$:
$$a(s_1)(\langle s_2, s_3 \rangle) = \langle [s_1, s_2]_A, s_3 \rangle + \langle s_2, [s_1, s_3]_A \rangle$$

4. what remains to make the flows to be automorphisms: they have to preserve the bracket $[,]_A$ (i.e. the very way the flows are produced by sections):
$$[s_1, [s_2, s_3]_A]_A = [[s_1, s_2]_A, s_3]_A + [s_2, [s_1, s_3]_A]_A$$

5. the mysterious one (the only hint now is that we cannot want $[s, s]_A = 0$ (i.e. $s$ to be conserved by its flow) if $\langle, \rangle$ is to be conserved):
$$[s, s]_A = \frac{1}{2} a^*(d\langle s, s \rangle)$$

Perhaps it is worth saying that all these properties are verified for (1) *without computation* (except for 5.). This is particularly convenient for the Jacobi identity 4.

## And how it can help us to find the global objects behind C. a.'s

The point is this. Suppose we did the same business for 1-dim. var. problems. Then instead of a C. a. we would arrive first at $A = \mathbb{R} \oplus TM$ and then we would realize that only the exact sequence $0 \to \mathbb{R} \to A \to TM \to 0$ and the Lie algebroid structure are natural. Starting with a 1-form on $M$ we would arrive at a connection on $A$, i.e. a splitting of the sequence. The global object is clear here – a principal $U(1)$ bundle. It is also clear from the physical point of view: a $U(1)$ bundle is exactly the thing that enables us to compute amplitudes for noncontractible closed curves.

This suggests that the global object behind an exact C. a. is a thing that allows computing amplitudes for general (i.e. noncontractible) closed surfaces. And this is known. There is even a book devoted to the answer by Brylinski (I guess you know it as the author thanks you warmly in the introduction).

A natural generalization comes to mind: why to use only 1 and 2 and not say 3? I really don't know what kind of geometry is hidden here (when compared e.g. with Poisson groupoids). The corresponding $A$ seems to be established easily, however. As a first step one uses $TM \oplus \bigwedge^2 T^*M$; the formula (1) remains (only $\theta$'s are now 2-forms) and $\langle, \rangle$ becomes a bilinear $T^*M$ valued form.

## Some trivialities on C. a.'s (classification of exact ones)

Notice that $aa^* = 0$ for any C. a. (because of 5. and 2.) and that Im $a^*$ can be characterized as the kernel of $\langle,\rangle$ restricted to the kernel of $a$. Also notice that $A/\text{Im } a^*$ (if $a^*$ has a constant rank) is a Lie algebroid.

If the sequence $0 \to T^*M \to A \to TM \to 0$ is exact, we call $A$ an exact algebroid. We'll see that locally there is only one such thing and globally they are classified by $H^3(M, \mathbb{R})$. To do it let us first notice that the sequence admit a global *isotropic* splitting $C \subset A$ (by a partition of unity, say); we shall call it a connection on $A$. For any maximally isotropic subbundle $C$ of $A$ (not necessarily a connection) we construct a section $c$ of $\bigwedge^3 C^*$ (the curvature of $C$) as follows: we take two sections of $C$ and observe that their bracket modulo $C$ is completely local, i.e. we have a map $C \otimes C \to A/C$. Then we observe that $A/C = C^*$ and check the complete skew-symmetry. $C$ is a Dirac structure iff $c = 0$.

If $C$ is a connection (i.e. transversal to $a^*(T^*M)$), then $c$ becomes a 3-form on $M$. Connections form an affine space over $\Gamma(\bigwedge^2 T^*M)$. If we add a 2-form $\beta$ to $C$, $c$ changes to $c + d\beta$. With this we notice that $c$ is closed (as $v(c) = d(v \lrcorner c)$) and thus can be locally annulated using an appropriate $\beta$. We see that $c$ gives us a class in $H^3(M, \mathbb{R})$ independent of $C$.

And this is really a classification: Algebroids with connection are classified by closed 3-forms and without connections by the 3rd de Rham cohomology. First notice that if we use $C$ to identify $A$ with $TM \oplus T^*M$, we have a formula like (1), except for an extra term $c(v_1, v_2, \cdot)$, so that $c$ really specifies the algebroid structure. One can also verify directly that for a closed $c$ the modified (1) defines a C. a. structure; if it is too tedious, one can use local annulations of $c$ and glue the trivial C. a.'s (with the formula (1)) using closed 2-form on overlaps.

## Why $[s, s]_A = \frac{1}{2} a^*(d\langle s, s \rangle)$?

Roughly speaking, because of a gauge transformation using the function $\frac{1}{2}\langle s, s \rangle$. Well, it needs some explanation.

The sections of the Lie algebroids associated to a principal $U(1)$ bundle $U \to M$ are the $U(1)$ invariant vector fields on $U$, i.e. the infinitesimal automorphisms of $U \to M$. We'd like to see the infinitesimal automorphisms of a Dixmier–Douady gerbe as the sections of a C. a. However, it is not that easy to define the infinitesimal automorphisms of a gerbe. For principal bundles it was simple, but if we defined bundles in the style of gerbes via their sheaves of sections (satisfying the obvious axioms), even here it would be difficult. *And suddenly, the explanation* Brylinski invented connective structure to *ends*.

## The BFG and Poisson–Lie T-duality

So it seems quite plausible that the global objects behind exact C. a.'s are the Dixmier–Douady gerbes, while the non-exact case remains unclear (a sort of IQ test: *U(1)-bundles : Lie groupoids = Dixmier–Douady gerbes : ?*). Now I'll try

to show that I'm far from real understanding the exact case while the distance between exact and non-exact is perhaps not that big.

The analogy with 1-dim case is again helpful. If we tried to define $U(1)$ bundles in the style of gerbes, we'd have the sheaves of sections or something of that kind. Of course, it would be equivalent, but we'd miss completely the fact that the bundle is first of all a *manifold*; compared with this, the $U(1)$ action is just an inessential ornament. This is nicely illustrated in mechanics. It seems to be largely accepted that the correct setting for mechanics (relativistic and prepared to quantization) is a principal $U(1)$ (or $\mathbb{R}$) bundle $U$ over a manifold $M$ (the extended phase space) with a connection such that the curvature $\omega$ has 1 characteristic direction everywhere. If one understands mechanics as the high-frequency limit of something, this bundle is completely natural. But still it is felt that it is better to work with $(M, \omega)$ and consider $U$ as a complication. But a lot of symmetry is missed in this way.

If we perceive $U$ as manifold, we see it has a $U(1)$ invariant distribution of hyperplanes with 1 characteristic direction and we want to find the characteristic curves. It may happen that we find a larger group of symmetry containing $U(1)$. In the $(M, \omega)$ approach one would not notice this symmetry unless $U(1)$ was in its center (and it is required to be a direct product to have a moment map). But in $U$ it does not matter. We can use the symmetry to reduce the problem (by contact reduction) and moreover (compare with symplectic reduction) solving the reduction we solve the complete problem (an application to symplectic reduction: to solve the problem we have to solve the reduction and to compute a single integral (to lift from $M_{red}$ to $U_{red}$) – this is the method of the reduced action).

And even more: we can use another $U(1)$ (or $\mathbb{R}$) in our group to make $U$ a principal bundle over another manifold. This gives us an equivalence of two mechanical models. But it is rather hard to imagine a sheaf over $M$ becoming a sheaf over another manifold.

The point is that this kind of equivalence exists in 2-dim problems as well. It is called Poisson–Lie T-duality. So it seems that there is some BFG (after Roald Dahl) and the usual definition of a gerbe is just its shadow on a manifold. Poisson–Lie T-duality enables us to throw shadows in different directions. This letter is becoming too long, but if you are interested in the duality, I can write you the details (and how general C. a.'s and Drinfeld's classification of Poisson homogeneous spaces enter).

<div style="text-align: right;">
Pavol Ševera<br>
IHÉS<br>
severa@ihes.fr
</div>

**Letter 2**

## "Nonabelian Noether theorem"

The starting point of PL T-duality was that 2d systems admit nonabelian conservation laws. It means that for each oriented curve on the surface we have an element of a Lie group $\tilde{G}$ (instead of $\mathbb{R}$) which is multiplicative with respect to concatenation and for contractible closed curves it is equal to 1. In other words, we have a $\tilde{\mathfrak{g}}$ valued 1-form satisfying Maurer–Cartan equation.

Now it would be nice to see an analogue of Noether theorem. Suppose again that the system is simply given by a 2-form $\alpha$ on $M$. To produce a nonab. c. law imagine a Lie group $G$ acting on $M$ and an identification of $\mathfrak{g}^*$ with $\tilde{\mathfrak{g}}$ (the actions in this text are right, if they are not evidently left). Now we mimic Noether theorem: we set $\beta = v \lrcorner \alpha$ (where $v \in \Gamma(\mathfrak{g}^* \otimes TM)$ is the infinitesimal action of $G$, so that $\beta$ is a $\mathfrak{g}^* = \tilde{\mathfrak{g}}$ valued 1-form). We would like $d\beta = \beta^2$ (the sign may be wrong here) on extremal surfaces and it is achieved if

$$v(\alpha) = \beta^2. \qquad (1)$$

Of course, it remains to understand the meaning of this equation. From now on we suppose that $G$ is a PL group and the identification $\mathfrak{g}^* = \tilde{\mathfrak{g}}$ is the standard one (it will become clear that this is more or less necessary). We let $G$ act on $A = TM \oplus T^*M$ in the following way: for a $g \in G$ we take the obvious action and we change it a bit. Namely, for any $x \in M$ we have a bivector $b$ at $x^g$ (this is the Poisson bivector at $g$ mapped by $G \to M$, $g \mapsto x^g$) and we use it to move $A_{x^g}$ in itself: a point $u + \theta \in A_{x^g}$ goes to $u + \langle b, \theta \rangle + \theta$. Now we see (1) as a statement of symmetry: the connection $C$ corresponding to $\alpha$ (i.e. the graph of $\alpha$) has to be preserved under this action. For example, if $\alpha$ is a symplectic form, (1) says that the action of $G$ is Poisson.

## The $M = G$ case

It is somewhat contradictory to take $M = G$ and at the same time the Lagrangian to be simply a 2-form on $M$ (because this was justified by the possibility to pass to a jet space). We shall do it before passing to general $M$'s, but it may be helpful to remark that any Lagrangian can be encoded in $A$ (by a lightlike hypersurface in the projective quadric given by the lightcone in $A_x$ – this is the general reparametrization-invariant case); if it is invariant with respect to the action of $G$ described above, we get a $\tilde{G}$ valued conservation law again.

The essential point is to use the Drinfeld double $D$; it explains nicely the strange $G$-action on $A$ and its connection with nonabelian conservation. (It is very simple and really natural: at the time we were inventing PL T-duality I did not know Poisson geometry at all, but it did not matter. The original idea was this: on extremal surfaces $\tilde{g} = P \exp \int \beta$ is a (locally) well-defined $\tilde{G}$ valued function and we can use it to lift the surface to $D$ via $g \mapsto g\tilde{g}$. We expected a duality between $G$ and $\tilde{G}$: the lifted surface might have been lifted from $\tilde{G}$ as well, using the supposed dual Lagrangian on $\tilde{G}$. The function $\tilde{g}$ is unique only up to multiplication by a constant from the right; the PDE for the lifted surfaces

1

must be right $\tilde{G}$-invariant. By the same argument they are right $G$-invariant, i.e. $D$-invariant. PL groups then fitted perfectly to this idea.)

We identify $M = G$ with $D/\tilde{G}$ (supposing it is possible) and $A = TM \oplus T^*M$ with $TD/\tilde{G} \cong T_G D$ ($TD/\tilde{G}$ will always mean $(TD)/\tilde{G}$; $TM \subset T_G D$ is obvious, $T^*M$ is vertical in $TD/\tilde{G}$ and the pairing with $TM$ is given by the inner product in $TD$). $G$ acts on $T_G D$ it the obvious way, and it is exactly the action on $A$ mentioned above. It means the following: if we encode $C$ as a $\tilde{G}$ invariant subbundle $C \subset TD$ (using $A = TD/\tilde{G}$), it will be $D$-invariant (as a special case, this contains the classification of Poisson-homogeneous symplectic structures). Surfaces $\Sigma \subset D$ satisfying $T\Sigma \subset C$ are exactly the lifts of extremal surfaces as mentioned above. We can exchange the roles of $G$ and $\tilde{G}$ and use the same $C$ to define a Lagrangian on $\tilde{G}$; extremal surfaces in $G$ can be exchanged with those in $\tilde{G}$ using the lift-and-project procedure. This is PL T-duality.

$D$ invariance is not directly visible in $G$. If we move a surface in $D$, its projection in $G$ remains extremal, but the transformation is complicated (it uses integrals to find the lift to $D$). It is important to stress that this is not an ordinary hidden symmetry. The lift of an extremal to $D$ is not unique (it is only up to $\tilde{G}$-action), and the result of the $D$-action is different for different lifts. This is very similar to the contact symmetries described in the previous note.

It is also interesting to understand PL T-duality as a phase-space transformation. We have to choose boundary conditions to build up a phase space. First suppose the case of closed strings (i.e. cylinders). The lift of a cylinder is not necessarily a cylinder, but we confine ourselves to cylinders in $D$. Doing this we are constrained to a submanifold of the total phase space and the symplectic form may become degenerate. And it actually becomes, making PL T-duality into a well-defined map on the reductions by null-directions. In this way, PL T-duality is a canonical transformation (in fact, the Poincaré–Cartan invariant can be easily expressed for the lifted strings and the subgroups $G, \tilde{G} \subset D$ do not enter this expression – this is the proof).

For open strings (strips) we have the constrain $u \lrcorner \alpha = 0$ for $u$ tangent to the boundary of the strip (this comes from the variational principle allowing variations of the boundary). As a result, the lifted strips have their boundaries on cosets of $G$ (of the type $Gd$, $d \in D$). The dual picture is a bit strange: the endpoints of the strings are constrained to move in the dressing orbits (D-branes) and they feel there an elmg force equal to the symplectic form on these orbits. We may also choose an isotropic subgroup $L \subset D$ (with the dimension equal to that of $G$) and consider strips with boundaries attached to the cosets of $L$. The strings in $G$ now feel the D-branes and their elmg fields given by the corresponding Dirac structure. In fact, this is a meaning of Dirac structures in 2d var. problems – they are D-branes.

If we choose two cosets of $L$ and consider PL T-duality only for the strips bounded to these two, it becomes a canonical transformation again. Finally, if we again consider open strings without D-branes, we can remove non-uniqueness of lifting by requiring one side of the lifted strip to lie in $G$. In this way the right action of $G$ on the phase space is well-defined and this action is Poisson.

## The $M = D/\tilde{G}$ case

The aim of this section is to explain its title. In the previous section we supposed $D$ to be diffeomorphic to $G \times \tilde{G}$ (or we may have worked locally). If we want the projections to be well-defined and to take the right $D$-symmetry seriously (and we should, as it appears to be at the heart of the thing), we are led to $M = D/\tilde{G}$. To find a variational problem on $M$, we'd like to identify $TM \oplus T^*M$ with $TD/\tilde{G}$, i.e. to find a splitting of $T_xD$. The picture shows the solution of the previous section:

(picture missing)

A problem appears if the intersection at $x$ is not transversal. If we use another isotropic subgroup $L$ instead of $G$, we get another splitting and thus another $\alpha$ on $M$. But these two splittings differ just by adding a closed 2-form. So this solves the problem: we divide $M$ into patches and for each we choose an appropriate $L$ (for example a conjugate of $G$ – this always works). In this way we find an exact CA structure on $TD/\tilde{G} \to M$ and the variational problem on $M$ can be formulated. The CA structure can be described as follows: for the right-$D$-invariant sections of $TD/\tilde{G}$ the bracket is just the Lie bracket on $\mathfrak{d}$; for arbitrary sections we use the Leibniz rule and the rule for the symmetric part of the bracket. The corresponding element of $H^3(M, \mathbb{R})$ gives us an obstruction for local transitivity of the action of $G$ on $D/\tilde{G}$. A closed 3-form on $M$ with this cohomology can be easily computed using a connection on $D \to M$.

An intriguing thing is that we do not need $G$ – it does not have to exist. The final definition of CA structure on $TD/\tilde{G}$ does not mention $G$ at all. At the very beginning an action of $G$ on $M$ was in the heart of the thing, but now we see it is not needed at all. The heart is the $D$-symmetry. I can't say I understand it completely, but the analogy with contact symmetries in mechanics is very usefull: there we did not need a subgroup transversal to $U(1)$, although we might have used one.

Now PL T-duality looks as follows: we have a group $D$ with invariant inner product and we choose a right-invariant $C$ (or a field of lihgtlike hypersurfaces mentioned above). For any isotropic half-dimensional (and closed) $\tilde{G} \subset D$ (this resembles quasi-Hopf algebras) we have a variational problem on $D/\tilde{G}$ and all these problems are connected by PL T-duality.

## The general case

There was a Courant algebroid hidden in the previous sections – it was $\mathfrak{d}$. Generally, we have a transitive Courant algebroid $A \to N$. We take the Lie algebroid $A/\operatorname{Im} a^*$ and assume that there is a corresponding principal bundle $D \to P \to N$. We pick up a $\tilde{G} \subset D$ and put $M = P/\tilde{G}$. This is more-or-less all, but it hides nice topology that I haven't explored as yet, except for simple cases.

**Towards classification of CA's**

Recall that exact CA's $A \to M$ (for which $0 \to T^*M \to A \to TM \to 0$ is exact) are classified by $H^3(M; \mathbb{R})$.

Local classification of $A$'s with constant-rank anchors is reduced to local classification of transitive $A$'s (living on the orbits of $A$) and their families. We shall see that transitive $A$'s are locally classified by Lie algebras with invariant inner products (nondegenerate symmetric bilinear forms).

If $A$ is an arbitrary CA and the anchor $a$ has constant rank, $A' = A/\operatorname{Im} a^*$ is a Lie algebroid and the kernels of $A' \to TM$ are Lie algebras with invariant inner product (the product is invariant with respect to flows in $A'$). Suppose that $A$ (and thus $A'$) is transitive and that $A'$ can be trivialized, i.e. $A' = TM \oplus \mathfrak{g}$. Then the preimage $A_0$ of $TM$ upon $A \to A'$ is an exact Courant subalgebroid of $A$. We may decompose $A$ orthogonally to $A_0 \oplus \mathfrak{g}$. Let $t$'s denote sections of $\mathfrak{g}$ (maps $M \to \mathfrak{g}$) and $s$'s sections of $A_0$. It is easy to see that

$$[s_1, s_2]_A = [s_1, s_2]_{A_0} \tag{1}$$

$$[s, t]_A = -[t, s]_A = \mathcal{L}_{a_0(s)} t \tag{2}$$

$$[t_1, t_2]_A = [t_1, t_2] + a_0^*(\langle dt_1, t_2 \rangle) \tag{3}$$

($a_0$ is the restriction of $a$ to $A_0$). And vice versa, if $A_0$ is an exact CA and $\mathfrak{g}$ a Lie algebra with an invariant inner product, $A_0 \oplus \mathfrak{g}$ becomes a CA in this way.

Since any transitive Lie algebroid can be locally trivialized and any exact CA is locally the standard $TM \oplus T^*M$, we have a local classification of transitive CA's.

Let us look at $A = TM \oplus T^*M \oplus \mathfrak{g}$ when $M$ is a circle. The sections of any CA form a Loday algebra and to get a Lie algebra we have to mod out the image of the symmetric part of the bracket, i.e. $a^*$ of exact 1-forms. In our case it means that we have to replace the 1-form part by its integral over the circle. From (3) we get

$$[t_1, t_2]' = [t_1, t_2] + \oint \langle dt_1, t_2 \rangle, \tag{4}$$

so that the standard central extension of loop algebras is born here. The Lie algebra of sections of $A$ modulo exact forms is the standard semidirect product of this central extension and of the algebra of vector fields on the circle. I wonder if this fact can be connected with the role of CA's in 2d variational problems.

It is quite easy to go beyond local classification in the case of transitive CA's. Perhaps you'd prefer to work it out yourself before reading the following incomplete and a bit confused discussion (and send me the result).

The point is that we can trivialize $A$ locally ($A = TM \oplus T^*M \oplus \mathfrak{g}$) and we have to look what happens when we glue these local trivializations. Recall the case of exact CA's. There we used closed forms on overlaps to glue the trivializations (the group (or sheaf of groups) of closed 2-forms acts by automorphisms on $TM \oplus T^*M$ (in fact, on any exact CA)). We needed a 1-cocycle in this sheaf; the corresponding cohomology is isomorphic to $H^3(M; \mathbb{R})$. It turned out that

this cohomology classifies exact CA's. A nice abstract nonsense description is this: there is 1-1 correspondence between exact CA's and $\mathcal{Z}^2$-torsors ($\mathcal{Z}^2$ is the sheaf of closed 2-forms); the correspondence is given via the action of $\mathcal{Z}^2$ on exact CA's.

For exact CA's we have $\mathfrak{g} = 0$; we have to understand what replaces the "glueing group" $\mathcal{Z}^2$ if $\mathfrak{g} \neq 0$. Let us describe its Lie algebra (sheaf of Lie algebras, more precisely). It is a central extension of the sheaf $\mathfrak{g}$ (the smooth $\mathfrak{g}$-valued functions) by $\mathcal{Z}^2$. It comes from (3) when "considered modulo closed 1-forms" (i.e. by taking differential of the 1-form on the right-hand side of (3)):

$$[t_1, t_2]_c = [t_1, t_2] + \langle dt_1 \stackrel{\wedge}{,} dt_2 \rangle. \tag{5}$$

This algebra appears for $A = A_0 \oplus \mathfrak{g}$ is the following way. Consider sections of the kernel of $A$, i.e. of $T^*M \oplus \mathfrak{g}$ and the corresponding flows (infinitesimal automorphisms) of $A$. As a result, we have an action of the central extension (5) by automorphisms on $A$ (the flow generated by a 1-form $\alpha$ depends only on $d\alpha$ – that's how 2-forms appear).

A bit imprecisely, transitive CA's are classified by torsors of these central extensions. A precise statement is desirable, of course, but I let it here as it is.

And the final remark. Suppose that $A$ is a transitive CA and that the Lie algebroid $A' = A/\operatorname{Im} a^*$ integrates to a principal bundle $G \to P \to M$. The bundle $\pi^* A$ (where $\pi : D \to M$ is the projection) seems to be quite interesting. Let us just remark that if $H \subset G$ is a closed maximally-isotropic subgroup and $\pi_H : P/H \to M$ is the projection then $\pi_H^* A$ is an exact CA over $P/H$. The CA structure appears as follows: for pullbacks of setions of $A$ to $\pi_H^* A$ we let the bracket to be that on $A$. This specifies the bracket of any sections of $\pi_H^* A$. This example is important in PL T-duality, that is why I mention it here.

**<span style="color:red">Letter 4</span>**   **Some topological nonsense**

If $G$ is a Lie group, $\underline{G}$ will denote the sheaf of smooth maps to $G$. If the Lie algebra $\mathfrak{g}$ is provided with an $Ad$-invariant inner product $\langle,\rangle$ (of any signature), we shall describe a sheaf of groups $\underline{\tilde{G}}$. It is a central extension of $\underline{G}$ by the sheaf of closed 2-forms $\mathcal{Z}^2$.

The expression $\langle u,[v,w]\rangle$ defines a both-sides invariant (and thus closed) 3-form $c$ on $G$. The sections of $\underline{\tilde{G}}$ are pairs $(g,\alpha)$, where $g$ is a smooth map to $G$ and $\alpha$ is a 2-form satisfying $d\alpha = g^*c$. The product of two sections is defined by
$$(g_1,\alpha_1)\cdot(g_2,\alpha_2) = (g_1g_2, \alpha_1+\alpha_2+\langle g_1^{-1}dg_1 \stackrel{\wedge}{,} dg_2\, g_2^{-1}\rangle).$$

Suppose we are given a principal $G$-bundle $P \to M$, i.e. an element of $H^1(M,\underline{G})$, and we try to extend it to a $\underline{\tilde{G}}$-torsor, i.e. an element of $H^1(M,\underline{\tilde{G}})$. The obstruction lies in $H^2(M,\mathcal{Z}^2) = H^4(M,\mathbb{R})$ and is equal to the "1st Pontryagin class" $p_1(P)$, i.e. to the class of the 4-form $\langle \Omega \stackrel{\wedge}{,} \Omega \rangle$, where $\Omega$ is the curvature of a connection on $P$.

The group $H^1(M,\mathcal{Z}^2)$ acts on $H^1(M,\underline{\tilde{G}})$; the orbits of this action are the fibres of the map $H^1(M,\underline{\tilde{G}}) \to H^1(M,\underline{G})$. It is desirable to know the isotropy groups of this action. Let us first describe the isotropy group of the trivial element of $H^1(M,\underline{\tilde{G}})$. It is the kernel $K$ of $H^1(M,\mathcal{Z}^2) \to H^1(M,\underline{\tilde{G}})$, i.e. the group consisting of the classes $g^*[c]$, where $[c] \in H^3(G,\mathbb{R})$ denotes the cohomology class of $c$ (we use $H^1(M,\mathcal{Z}^2) = H^3(M,\mathbb{R})$) and $g$ runs through maps $M \to G$ (or through their homotopy classes). Thus, the extensions of the trivial $G$-bundle $G \times M$ form the group $H^3(M,\mathbb{R})/K$. For example, if $M$ is a closed oriented 3-fold and $G = SU(2)$, then $H^3(M,\mathbb{R})/K \cong \mathbb{R}/Z$. Moreover, any $P$ is trivial in this case, so that $H^1(M,\underline{\tilde{G}}) \cong \mathbb{R}/Z$.

The isotropy group for a general $P$ is similar. Again, it is the image of a homomorphism from the gauge group of $P$ to $H^3(M,\mathbb{R})$ (it is a twisted connecting morphism). It admits a purely differential description (based on lengthy computation I don't trust completely; the straightforward definition is obvious, but it uses e.g. Čech cocycles). Choose a connection on $P$ and let $\tau$ be the corresponding Chern–Simons 3-form. For any gauge transformation of $P$ take the horizontal part of $\tau$ with respect to the transformed connection. It gives us a closed 3-form on $M$ and the promised homomorphism.

## What the nonsense is good for (classification of transitive CA's)

If there were an equivalence between principal bundles and transitive Lie algebroids, the previous section would be about classification of transitive CA's. It is not the case, but actual classification of transitive CA's is quite close. If $A$ is a transitive CA, $A' = A/\operatorname{Im} a^*$ is a transitive Lie algebroid with invariant inner product on the vertical algebras. If $A'$ is given in advance, an $A$ exists iff $p_1(A') = 0$. Finally, the group $H^3(M,\mathbb{R})$ acts transitively on their types with the isotropy group given as before, but with the gauge group replaced by the group of automorphisms of $A'$ preserving the inner product.

1

The classification is derived by a simple (though a bit dull) computation. I don't know any more conceptual proof. However, if we insist on principal bundles, there is such a proof, decribed in the following section. It is natural to expect that if there is a global object behind $A$, $A'$ should integrate to a principal bundle.

For completeness, here are some details of the computation. We choose a connection on $A'$, so that we can write $A' = V \oplus TM$, where $V$ is the bundle of vertical Lie algebras. For a while, let $A_0$ denote the preimage of $TM$ under $A \to A'$. We identify $V$ with $A_0{}^\perp$. We split the exact sequence $0 \to T^*M \to A_0 \to TM \to 0$ by an isotropic subbundle. As a result, we have $A = T^*M \oplus TM \oplus V = T^*M \oplus A'$, with $TM$ and $T^*M$ isotropic and $V = (TM \oplus T^*M)^\perp$. Now we simply write down the general form of the Courant bracket $[,]_A$; $\alpha$ denotes a 1-form, $v$ a vector field, $t$ a section of $V$; $\nabla$ is the covariant derivative, $\Omega$ the curvature and $c$ is some 3-form to be explained below.

$$[t, \alpha]_A = [\alpha, t]_A = 0$$
$$[v, t]_A = -[t, v]_A = [v, t]_{A'} - \langle v \lrcorner \Omega, t\rangle$$
$$[t_1, t_2]_A = [t_1, t_2]_{A'} + \langle \nabla t_1, t_2 \rangle$$
$$[v, \alpha]_A = v(\alpha)$$
$$[\alpha, v]_A = -v \lrcorner d\alpha$$
$$[\alpha_1, \alpha_2]_A = 0$$
$$[v_1, v_2]_A = [v_1, v_2]_{A'} + c(v_1, v_2, \cdot)$$

When considered modulo $T^*M$, this $[,]_A$ gives $[,]_{A'}$, as it should. The $T^*M$-part of the right-hand sides then follows immediately from the axioms of CA's. The only freedom is in the last formula, where $c$ is some 3-form.

These fomulas give a CA provided $dc = \langle \Omega \stackrel{\curlywedge}{,} \Omega \rangle$. This comes from Jacobi identity for $[,]_A$, and up to now, it is the only mindless (and fortunately, quite short) computation.

Finally, we have to find out when two $c$'s give isomorphic $A$'s. It is certainly the case if $c_1 = c_2 + d\beta$ for some 2-form $\beta$: recall the arbitrary isotropic splitting of $0 \to T^*M \to A_0 \to TM \to 0$ that may be changed by an arbitrary 2-form. Now this arbitrariness is taken into account, and for a fixed $A'$ we have classification of $A \to A'$ (i.e. the isomorphisms between $A$'s are supposed to form commutative triangles with $\to A'$). However, an actual classification of $A$'s has to go a bit further: we are only given the isomorphism class of $A'$, and so we have to take into account the group of automorphisms of $A'$. The result has already been mentioned (this is the computation I wouldn't trust).

## Reductions

Recall that $[s, s]_A = a^* d\langle s, s\rangle/2$. If we want to represent a Lie algebra $\mathfrak{h}$ by sections of $A$ then these sections must have constant inner products. In this way $\mathfrak{h}$ becomes endowed with a possibly degenerate inner product.

Here is a particularly useful example. Let $G$ be as above (a Lie group with an invariant inner product). We take $A = TG \oplus T^*G$ twisted with the 3-form $c$ (i.e. we add $c : \bigwedge^2 TG \to T^*G$ to the standard Courant bracket). For any $v \in \mathfrak{g}$ we take $v + \langle v, \cdot \rangle$ and translate it right-invariantly; in this way we obtain a representation of $\mathfrak{g}$ in $\Gamma(A)$. For an alternative description (to be important later), take $A = TG \oplus TG$, where the first $TG$ is taken with $\langle,\rangle$ and the second one with $-\langle,\rangle$. The bracket of right-invariant fields of the first copy is defined to be simply the commutator (again in the first $TG$), and the same for left-invariant fields of the second $TG$. This bracket is uniquely extended to a Courant bracket on $A$.

Let a Lie algebra $\mathfrak{h}$ with a possibly degenerate inner product be represented by sections of $A$ via map $\phi : M \times \mathfrak{h} \to A$. Applying $a$ on these sections we have an integrable singular distribution on $M$ (and an action of the local group $H$, if $\mathfrak{h}$ is finite-dimensional). Suppose that $\phi$ has constant rank and the the space of leaves of this distribution is a manifold $N$. Then the $\mathfrak{h}$-invariant sections of $(\operatorname{Im} \phi)^\perp / \operatorname{Im} \phi \cap (\operatorname{Im} \phi)^\perp$ are sections of a CA $A_{red} \to N$.

As a special case, let $\pi : P \to M$ be a principal $G$-bundle, $\tilde{A} \to P$ an exact CA and let the usual representation of $\mathfrak{g}$ by vector fields be lifted to $A$, inverting the sign of the inner product. Then $\tilde{A}_{red} \to M$ is a transitive CA and the Lie algebroid $(\tilde{A}_{red})'$ is canonically isomorphic to the Atiyah sequence of $P$. And vice versa, for any $A \to M$ with $A'$ isomorphic to $TP/G$ there is a canonical $\tilde{A}$ with $\tilde{A}_{red} = A$. This seems to be important. Supposing we understand the global objects behind exact CA (a seemingly tractable problem), we should understand them for transitive CA's as well.

Here is a sketch of construction of $\tilde{A}$ for a given $A$, $A' = TP/G$. It is a generalization of $TG \oplus TG$ construction mentioned above (in that case $M$ was a point). Put $\tilde{A} = \pi^* A \oplus VP$, where $VP = P \times \mathfrak{g}$ is the vertical part of $TP$. The anchor of $\tilde{A}$ is obvious (recall $A' = TP/G$). The inner product on $VP$ is $-\langle,\rangle$ (the inner product on $\mathfrak{g}$), on $\pi^* A$ it is as on $A$ and $VP \perp \pi^* A$. Finally the Courant bracket: for constant sections of $VP$ (i.e. for elements of $\mathfrak{g}$) it is the commutator in $\mathfrak{g}$ and for pullbacks of sections of $A$ it is the pullback of their bracket. This information specifies $[,]_{\tilde{A}}$ uniquely.

Now we can return to the first section. As it turns out, specifying an exact $\tilde{A}$ over $P$ with a representation of $\mathfrak{g}$ in $\tilde{A}$ as above is equivalent to extending $P$ to a $\underline{\tilde{G}}$-torsor. It is easily seen: locally, for nice $U \subset M$ and a trivialization $P|_U = G \times U$, $\tilde{A}$ can be trivialized to $(TG \oplus TG) \oplus (TU \oplus T^*U)$; glueing on overlaps is given by sections of $\underline{\tilde{G}}$, extending the glueing of $P$ by sections of $\underline{G}$.

And finally, reductions seem to be the proper framework for PL T-duality. If we take a half-dimensional isotropic subgroup $H \subset G$, the reduction of $\tilde{A}$ with respect to $H$ gives an exact CA $A_H \to P/H$, where variational problems can live.

## A fairy-tale

I have to admit that whatever I can firmly say about these topics is of rather trivial nature. It is the possible connection *CA's and their global objects (gerbes?)*

– *quantum groups* – *Chern–Simons theory* – *quantum PL T-duality* – *etc.* that promises something deeper. There may be some new light in this note, together with some new confusion.

# Courant algebroids and Dixmier–Douady gerbes

**Abstract** Exact Courant algebroids emerge as infinitesimal symmetries of Dixmier–Douady gerbes with connective structure. Transitive CA's are also accessible.

## Exact CA's as conducting bundles

We make the following notation: if $W$ is a vector space then $WM$ denotes the trivial vector bundle $W \times M$ over $M$.

Let $0 \to \mathbb{R}M \to B \to TM \to 0$ be an exact sequence of vector bundles and denote the projection $B \to TM$ as $\pi_B$. If it is an Atiyah sequence then every section $t$ of $B$ would generate a flow of $B$, covering the flow of $TM$ generated by $\pi_B(t)$. In particular, a section $f$ of $\mathbb{R}M$ would let $b \in B$ flow with velocity $\dot{b} = df(\pi_B(b)) \in \mathbb{R}$.

However, if $0 \to \mathbb{R}M \to B \to TM \to 0$ is just an exact sequence of vector bundles, there is no natural way of generating flows by sections of $B$. The sequence admits a lot of automorphisms (automorphisms of $B$ preserving the sequence). There is 1-1 correspondence between these automorphisms and 1-forms on $M$: if $\theta$ is a 1-form then $b \in B$ gets mapped to $b + \theta(\pi_B(b))$. We could compose a tentative flow of $B$ generated by some section $t$ with any such an automorphism, at least if $t$ itself is conserved, i.e. if $\theta(\pi_B(t)) = 0$.

As we see, sections of $B$ do not know how to generate flows of $B$. It is natural to invent the concept of a *conducting bundle*, whose sections would explain sections of $B$ how to do it, in all possible ways. By definition, a conducting bundle for $B$ is an exact sequence $0 \to T^*M \to A \to TM \to 0$ (with $A \to TM$ denoted as $\pi_A$). If $t$ is a section of $B$ then its *conductor* is a section $s$ of $A$ such that $\pi_A(s) = \pi_B(t)$; if we choose a conductor then $t$ (being instructed by the conductor) generates a flow of $B$ and if we change the conductor by a 1-form $\theta$ then the flow gets changed by $\dot{b} = \theta(\pi_B(b))$. Finally, if $f$ is a function (a section of $\mathbb{R}M$) then we want the flow generated by $f$, with 0 as the conductor, to be $\dot{b} = df(\pi_B(b))$.

It is good to restate the definition in a bit more formal way. Temporarily, let $C$ denote the kernel of $\pi_A - \pi_B$ contained in $A \oplus B$. We have the obvious map $\pi_C : C \to TM$. By a flow of $B$ covering a vector field $v$ we mean a vector field on $B$ such that 1. its local flow maps linearly fibres to fibres, 2. it can be projected via $\pi_B$ to $TM$ yielding the vector field in $TM$ generated by $v$ and 3. it can be restricted to $\mathbb{R}M \subset B$, yielding the natural lift of $v$ to $\mathbb{R}M$. What we want is a linear map $\rho$ from sections of $C$ to flows of $B$ such that 1. if $u$ is such a section then $\rho(u)$ covers $\pi_C(u)$, 2. if $\theta$ is a 1-form, considered as a section of $C \subset A \oplus B$, then $\rho(\theta)$ is given by $\dot{b} = \theta(\pi_B(b))$ and similarly, 3. if $f$ is a function then $\rho(f)$ is given by $\dot{b} = df(\pi_B(b))$. Notice the possibility of transfering $f \in \Gamma(B)$ to $df \in \Gamma(A)$.

It is more or less clear from the definition that we can speak about *the* conducting bundle of $B$ instead of *a* conducting bundle: the sections of $A$ are exactly all the ways of conducting sections of $B$ in their attempts to generate flows (a simple construction of $A$ will be given bellow). Therefore, the assign-

1

ment $B \mapsto A$ is functorial. We can ask how automorpisms of $B$ (1-forms) are mapped to automorphisms of $A$. As one easily sees, a 1-form $\theta$ is mapped to $s \mapsto s + \pi_A(s) \lrcorner d\theta$, $s \in A$. In particular, if $\theta$ is closed, $A$ remains intact.

As we shall see now, $A$ is in fact a Courant algebroid. First, we shall construct a bracket on sections of $A$, i.e. we shall show how a section $s$ of $A$ generates a flow of $A$ covering $\pi_A(s)$. It will come from a flow of $B$ by functoriality. Although $s$ alone does not give a flow of $B$, since we also need to specify a $t$ with $\pi_B(t) = \pi_A(s)$, any two such $t$'s differ just by a function. And the flow generated by a function $f$ yields zero flow of $A$ (as stated above, with $\theta = df$ being closed).

Now we shall construct an inner product on $A$ such that if $\theta$ is a 1-form and $s$ a general section of $A$ then $\langle \theta, s \rangle = \theta(\pi_A(s))$ (actually, it will reveal the meaning of the inner product and of the formula $[s,s]_A = d\langle s, s \rangle/2$). Take a section $s$ of $A$ and a section $t$ of $B$ with $\pi_B(t) = \pi_A(s)$ and consider the corresponding flow of $B$. We can ask what is the rate of change of $t$ under this flow. Generally, we cannot expect it to be zero, but at least we know that it is a section of $\mathbb{R}M$. If we add a 1-form $\theta$ to $s$, this rate gets changed by $\theta(\pi_B(t)) = \theta(\pi_A(s))$. With a little effort one shows that the rate is independent of $t$ (i.e. it depends only on $s$) and in fact is given by a quadratic form on $A$ as $\langle s, s \rangle/2$.

Now we see that $A$ is a Courant algebroid. The flows of $A$ generated by its sections have to preserve all its structure, so the only thing to be checked is $[s,s]_A = d\langle s, s \rangle/2$. Just recall the definition of $\langle , \rangle$: since $t$ is changing with velocity $\langle s, s \rangle/2$, $s$ must be changing with velocity $d\langle s, s \rangle/2$.

We shall discuss few additional topics. An Atiyah-sequence structure on $B$ is the same as a flat connection on $A$ (flat connection is a splitting of $0 \to T^*M \to A \to TM \to 0$ by a Dirac structure) – we take horizontal sections of $A$ as conductors. If the connection is not flat, its curvature (a closed 3-form) will appear as the Jacobiator of the bracket on $B$.

If we take $B = \mathbb{R}M \oplus TM$, we can easily construct its conducting bundle, setting $A = T^*M \oplus TM$. $B$ carries standard (trivial) Atiyah-sequence structure. The sections of $TM \subset A$ will conduct this structure; this specifies the conducting-bundle structure on $A$ uniquely (since every section of $A$ is a section of $TM$ plus a 1-form). This allows us to construct conducting bundle for any $B$ in a functorial way: any $0 \to \mathbb{R}M \to B \to TM \to 0$ admits a splitting (i.e. is isomorphic to $\mathbb{R}M \oplus TM$); any two splittings differ by a 1-form $\theta$ and we use $d\theta$ to change the splitting of $A$.

Let $L_{1,2}$ be Hermitean line bundles and $B_{1,2}$ their Atiyah bundles. The Atiyah bundle of their product $L_1 \otimes L_2$ will be denoted as $B_1 \vee B_2$. It can be defined as follows: we take the kernel of $\pi_{B_1} - \pi_{B_2}$ in $B_1 \oplus B_2$ and mod out $\mathbb{R}$ embedded antidiagonally into $\mathbb{R} \oplus \mathbb{R}$. This definition makes sense for any exact sequences $0 \to \mathbb{R}M \to B \to TM \to 0$ (not only for Atiyah sequences). We may therefore ask how to construct $A_{B_1 \vee B_2}$ out of $A_{B_1}$ and $A_{B_2}$. It is very similar: we take the kernel of $\pi_{A_{B_1}} - \pi_{A_{B_2}}$ and mod out the antidiagonal $T^*M$. The result will be denoted as $A_{B_1} \vee A_{B_2}$. This definition makes sense for arbitrary exact CA's and makes them into a symmetric monoidal category (with unit the standard CA $T^*M \oplus TM$). If we pass to isomorphism classes, we end up with

the group $H^3(M, \mathbb{R})$.

Finally, if $0 \to \mathbb{R}M \to B' \to TM \to 0$ is an Atiyah sequence then $A_{B'}$ is canonically isomorphic to $T^*M \oplus TM$, since the sections of $B'$ need no conduct. And if $0 \to \mathbb{R}M \to B \to TM \to 0$ is just an exact sequence of vector bundles then $A_{B \vee B'}$ is canonically isomorphic to $A_B$.

## Where nontrivial CA's come from

CA's appearing in the previous section are trivial (isomorphic to $T^*M \oplus TM$), since every sequence $0 \to \mathbb{R}M \to B \to TM \to 0$ admits a splitting. Yet, it reveals the meaning of these CA's. The way out has already been mentioned: $A_{B \vee B'}$ is canonically isomorphic to $A_B$ if $B'$ is an Atiyah sequence. Or alternatively (and more or less equivalently), an automorphism of $B$ defined by a closed 1-form yields identity in $A_B$. In other words, specifying $A$ we do not specify $B$ completely and it may result in nontrivial global behaviour.

Imagine $M$ is covered by open subsets $U_i$ and on each $U_i$ we have a sequence $0 \to \mathbb{R}U_i \to B_i \to TU_i \to 0$. On each overlap $U_i \cap U_j$ we choose an isomorphism between $B_i$ and $B_j$. We do not suppose that these isomorphisms form a cocycle (yielding a global $B$ glued from $B_i$'s) but we require its coboundary to consist of *closed* 1-forms. Although $B_i$'s do not glue into a global bundle, $A_{B_i}$'s do.

Alternatively, suppose we have Atiyah sequences $B_{ij}$ on every $U_i \cap U_j$ and isomorphisms $B_j \cong B_i \vee B_{ij}$. $B_{ij}$'s should form a cocycle in the following sense: for any $ijk$ we choose an isomorphism between $B_{ij} \vee B_{jk} \vee B_{ki}$ and the trivial Atiyah sequence on $U_i \cap U_j \cap U_k$ and these isomorphisms should be coherent for any $ijkl$. Again, $A_{B_i}$'s glue to a global $A$. To come back to the previous picture, suppose that all $U_i \cap U_j$'s are contractible and choose an arbitrary trivialization of $B_{ij}$.

## Dixmier–Douady gerbes

They are defined in Brylinski's book as gerbes with the band $\underline{U(1)}$, where $\underline{U(1)}$ is the sheaf of smooth maps to $U(1)$. I briefly recall the relevant definitions. A *presheaf of groupoids* has obvious meaning (restrictions are supposed to compose only up to natural transformations, but this subtlety plays no role in the sequel). If $\mathcal{C}$ is a presheaf of groupoids and $a, b \in \mathcal{C}(U)$ then morphisms between restrictions of $a$ and $b$ organize to a presheaf on $U$ denoted $\text{Hom}(a, b)$. $\mathcal{C}$ is a *stack* (or a sheaf of groupoids) if 1. every $\text{Hom}(a, b)$ is a sheaf and 2. for any covering $\{U_i\}$ of $U$, whenever we choose objects $a_i \in \mathcal{C}(U_i)$ and isomorphisms $\phi_{ij}$ beween restrictions of $a_i$ and $a_j$ to $U_i \cap U_j$ in a coherent way, there exist an object $a \in \mathcal{C}(U)$, unique up to canonical isomorphism, and isomorphisms $\phi_i : a|_{U_i} \to a_i$, commuting with $\phi_{ij}$. A *gerbe* is a stack where $\mathcal{C}(U_i)$ is not empty for some (fine enough) covering of $M$ and any two objects are locally isomorphic. $\mathcal{C}$ is a *Dixmier–Douady gerbe* (DDG) if isomorphism groups of its objects are isomorphic to $\underline{U(1)}$ in a coherent way. Finally, a *connective structure* on a DDG $\mathcal{C}$ is a functor assigning to $a \in \mathcal{C}(U)$ an exact sequence of vector bundles $0 \to \mathbb{R}U \to B(a) \to TU \to 0$ in such a way that an isomorphism of $a$, i.e. a map

3

$h$ to $U(1)$, is mapped to the 1-form $-\sqrt{-1}h^{-1}dh$ (viewed as an automorphism of $B(a)$).

The basic (or trivial) example of a DDG is the groupoid of principal $U(1)$ bundles (with usual isomorphisms). It has a canonical connective structure, producing Atiyah sequence out of a $U(1)$ bundle. Actually, every DDG with a global object is isomorphic to the trivial DDG. Maybe the best view of a DDG is a kind of 2-$U(1)$-bundle, with the abelian gauge group $U(1)$ replaced by the symmetric monoidal category (or more precisely, symmetric groupal groupoid) of principal $U(1)$-bundles with their obvious multiplication.

If we have a DDG with a connective structure then $A_{B(a)}$'s glue to a global $A$. A *curving* of the DDG is a connection on $A$ (it was defined by Brylinski in a bit different way). A DDG with a connective structure is a global framework for quantization of 2dim field theories and curving plays the role of Lagrangian. On the other hand, exact CA's form a framework for 2dim variational problems (their sections play the role of infinitesimal symmetries) and again, a connection on a CA gives the Lagrangian. The picture is quite consistent, but a lot of things remain to be understood.

## Transitive CA's

All CA's were exact up to now. Let us recall the relation between exact and transitive CA's. Let $P \to M$ be a principal $G$-bundle and let $\mathfrak{g}$ admit an invariant inner product. Also, let $A_P$ be an exact CA over $P$ and let the representation of $\mathfrak{g}$ by vector fields in $P$ be lifted to a representation by sections of $A_P$, preserving the inner product. We can produce a transitive CA $A_M$ over $M$ by reduction: the sections of $A_M$ are $\mathfrak{g}$-invariant sections of $A_P$ orthogonal to $\mathfrak{g}$. And vice versa, if $A_M$ is a transitive CA and if the Lie algebroid $A_M/T^*M$ integrates to a principal bundle $P$, then $A_P$ and the representation of $\mathfrak{g}$ is uniquely specified.

We have to understand what a representation of $\mathfrak{g}$ by sections of $A$ means at the level of DDG's. Atually, the best thing would be to find directly the reduction of the global object of $A$ (the "gauge groupoid of the DDG"), but I do not know how to do it. Instead, let us just notice the following thing. DDG's with connective structure are used to construct $U(1)$ bundles over loop spaces. Suppose we have such a DDG over a manifod $M$, $A$ the corresponding exact CA and let $\hat{LM} \to LM$ be the $U(1)$ bundle. Suppose that $G$ acts on $M$ and therefore $LG$ acts on $LM$. We want to find a lift to an action of $\hat{LG}$ (the standard central extension) on $\hat{LM}$. As it turns out, an infinitesimal part of this lift is the representation of $\mathfrak{g}$ by sections of $A$.

I do not know what is the meaning of this fact for Poisson groupoids. Let us consider the simplest example – the double of a Lie bialgebra. We shall ignore topological problems. Let $\mathfrak{g}$ be the double itself and $\mathfrak{h}, \tilde{\mathfrak{h}} \subset \mathfrak{g}$ the dual pair. Of course, $\mathfrak{g}$ is a transitive CA over a point; the principal bundle over the point is simply $G$. In this case the DDG over $G$ is the one used by Brylinski to construct $\hat{LG}$, i.e. we end up with the action of $\hat{LG}$ on itself. This picture does not seem to be simple, but actually, there is an interrelation between loop groups an PL

groups, a kind of classical reason for appearance of quantum group in WZW models (Gawedzki). Over $LH$ the central extension $\hat{LG}$ is canonically trivial, therefore we have an inclusion $LH \subset \hat{LG}$ (and similarly, $L\tilde{H} \subset \hat{LG}$). Take two loops $\gamma_1 \in LH$, $\tilde{\gamma}_1 \in L\tilde{H}$ and their product $\gamma_1 \tilde{\gamma}_1$ in $\hat{LG}$. We may decompose $\gamma_1 \tilde{\gamma}_1$ as $\tilde{\gamma}_2 \gamma_2$ times an element $c \in U(1)$. As it turns out, $c = \exp i \int_\Sigma \omega$, where $\Sigma$ is a disk in $G$ with boundary $\gamma_1 \tilde{\gamma}_1$ and $\omega$ is the symplectic form of the double symplectic groupoid.

**Letter 6**    **From associative algebroids to associative algebras**

Perhaps it is enough to say

$$\text{assoc. algebras : assoc. algebroids = groups : groupoids.}$$

By definition, an assoc. algebroid is a vector bundle $\mathcal{A} \to M \times M$, together with maps $\mathcal{A}_{(x,y)} \otimes \mathcal{A}_{(y,z)} \to \mathcal{A}_{(x,z)}$ depending smoothly on $x, y, z \in M$ and associative in the obvious sense (one can also use a Lie groupoid in place of $M \times M$).

There is an associative product on $\Gamma(\mathcal{A} \otimes |\det T^*|^{1/2}(M \times M))$ defined by

$$\alpha * \beta(x, z) = \int_{y \in M} \alpha(x, y) \cdot \beta(y, z)$$

(we should make some restriction on behaviour of $\alpha$ and $\beta$ at infinity to make the integral convergent). This formula defines associative algebra stucture on the space of generalized sections as well; it contains the algebra of $\mathcal{A}$-pseudodifferetial operators (Lagrangian distributions with respect to the diagonal of $M \times M$) and the algebra $\mathcal{D}(\mathcal{A})$ of $\mathcal{A}$-differential operators (pseudodif. oper. with support (not just singular support) in the diagonal).

## Clifford algebroids, symbol calculus, and nilpotent Diracs

Let $V$ be a vector space with inner product. One constructs $Spin(V) \subset Cl(V)$ in the usual way. There is a natural $\mathbb{Z}/(2)$-graded involution $t: Cl(V) \otimes \mathbb{C} \to Cl(V) \otimes \mathbb{C}$ (i.e. $t^2 =$ the parity operator, $\lambda^t = \lambda$ for $\lambda \in \mathbb{C}$, and $(ab)^t = (-1)^{|ab|} b^t a^t$) defined uniquelly by $v^t = iv$ for $v \in V$. Clearly $g^t = g^{-1}$ for $g \in Spin(V)$. $Cl(V)$ is filtered, with the associated graded Poisson algebra $\bigwedge^* V$; there (i.e. if we pass to symbols) $t$ becomes $i^{\text{degree}}$.

As a generalization, let $A \to M$ be a vector bundle with inner product; we shall define its *Clifford algebroid* $Cl(A) \to M \times M$. Let $V$ be a vector space with inner product, isomorphic to the fibres of $A$, and let $P$ be the corresponding principal $SO(V)$-bundle (i.e. the points of $P$ are isomorphisms of $V$ with the fibres of $A$). Suppose that $P$ admits a lift to principal $Spin(V)$-bundle $\tilde{P}$. Then

$$Cl(A) = \tilde{P} \times \tilde{P} \times Cl(V)/Spin(V) \times Spin(V),$$

where the $Spin(V) \times Spin(V)$-action is $(p_1, p_2, a)^{(g_1, g_2)} = (p_1 g_1, p_2 g_2, g_1^{-1} a g_2)$ (recall $Spin(V) \subset Cl(V)$, so that $g_1^{-1} a g_2$ makes sense). Multiplication in $Cl(A)$ is given by $(p_1, p_2, a)(p_2, p_3, b) = (p_1, p_3, ab)$ (easily seen to be well defined).

When restricted to the diagonal of $M \times M$, $Cl(A)$ becomes the bundle of Clifford algebras of $A$. Although $Cl(A)$ depends on $\tilde{P}$, it can be defined over a neighbourhood of the diagonal without use of $\tilde{P}$: simply replace $\tilde{P} \times \tilde{P}$ in its definition by a neighbourhood of diagonal of $P \times P$ and put similar condition on $Spin(V) \times Spin(V)$. Certainly, $\mathcal{D}(Cl(A))$ does not need any spin structure.

We define a map $t: Cl(A) \otimes \mathbb{C} \to Cl(A) \otimes \mathbb{C}$ covering the flip map of $M \times M$ by $(p_1, p_2, a)^t = (p_2, p_1, a^t)$. It generates a $\mathbb{Z}/(2)$-graded involution

1

on the algebra of $\mathit{Cl}(A)$-valued generalized half-densities and its subalgebras (provided we choose the condition at infinity in a flip-invariant way).

Finally, there is a reasonable symbol calculus on $\mathcal{D}(\mathit{Cl}(A))$, i.e. there is a filtration on $\mathcal{D}(\mathit{Cl}(A))$ such that the associated graded Poisson algebra is graded-commutative. Moreover, the map $t$ can be used to define subprincipal symbols. Namely, using trivialization and coordinates to make it a differential operator on $\mathbb{R}^n$ with coefficients in $\mathit{Cl}(V)$, the order of a $\mathit{Cl}(A)$-diff. operator is the order of the coefficient (in $\mathit{Cl}(V)$) plus *twice* the order of the derivative; this definition is clearly independent of the choices.

These symbols are functions on the graded symplectic supermanifold $(\mathcal{E}, \omega)$ with $\deg \omega = 2$, constructed by Alan. I recall its construction. Take $T^*P \times \Pi V^*$ with $T^*P$ graded by even numbers and $\Pi V$ with the usual grading (the standard notation is $T^*[2]P \times V^*[1]$); $\mathcal{E}$ is its symplectic reduction at zero. Using a trivialization, $\mathcal{E}$ is locally $T^*M \times \Pi V^*$; this can be used to define the symbol as the symbol of the coefficient times the symbol of the derivative; again (quite clearly), this is independent of the choices.

When we pass to symbols, $t$ becomes $i^{\mathrm{degree}}$. An operator $D$ of degree $k$ will be called *selfadjoint* if $D^t = i^k D$. Any $D$ of degree $k$ and parity $k$ mod 2 is uniquely split as $D = D_1 + D_2$, where $D_{1,2}$ are selfadjoint and $\deg D_1 = k$, $\deg D_2 = k - 2$ (namely $D_1 = (D + i^{-k}D^t)/2$, $D_2 = (D - i^{-k}D^t)/2$). The symbol of $D_2$ is the *subprincipal symbol* of $D$.

Here is an application to Courant algebroids and nilpotent Dirac operators. Suppose $\theta$ is a cubic function on $\mathcal{E}$ satisfying $\{\theta, \theta\} = 0$, i.e. we are given a Courant algebroid structure on $A$. Let $D$ be the selfadjoint 3rd degree operator with symbol $\theta$ (because $D$ is of degree 3, the selfadjointness condition specifies it uniquely). Since $\{\theta, \theta\} = 0$, $[D, D]$ is of degree at most 2. We have $[D, D]^t = -[D^t, D^t] = [D, D]$. Hence $[D, D]$ is of degree 0, i.e. $D^2 = $ function. In other words, for any Courant algebroid, there is a canonical choice for generating nilpotent Dirac operator.

## Symbols via tangent groupoid

The aim of this section is to keep finite number of dimensions to the very last moment. We shall need assoc. algebroids over groupoids (not just over $M \times M$); this generalization is clear (just half-densities on $M \times M$ are replaced by half densities on $\alpha$-fibres times half-densities on $\beta$-fibres). We'll describe an algebroid $\mathcal{T}\mathit{Cl}(A)$ over the tangent groupoid $\tau M$; the symbol calculus described above is contained in $\mathcal{T}\mathit{Cl}(A)$.

I first recall the structure of $\tau M$. Its base is $M \times \mathbb{R}$; points in this $\mathbb{R}$ will be denoted as $\epsilon$. Morphisms of $\tau M$ never change $\epsilon$, i.e. $\tau M$ can be viewed as a 1-parameter family of groupoids. Over $\epsilon \neq 0$ the groupoid is just $M \times M$, while over 0 it is $TM$. There is an $\mathbb{R}^*$-action on $\tau M$ (scale transformations): for $\epsilon \neq 0$, $\epsilon$ is just mapped to $\epsilon/\lambda$, while over $\epsilon = 0$, $TM$ is multiplied by $\lambda$ ($\lambda \in \mathbb{R}^*$). There is a natural manifold structure on $\tau M$ for which this action is smooth.

2

Now we can describe $\mathcal{T}Cl(A) \to \tau M$. Over $\epsilon \neq 0$ it is $Cl(A)$; over 0, it is a new algebroid $Gr(A) \to TM$:

$$Gr(A) = (TP \times \bigwedge V)/TSO(V).$$

Here $TSO(V)$ (the tangent bundle of $SO(V)$) is a Lie group in the usual way (it is a semidirect product of $SO(V)$ with its Lie algebra). Its action on $TP$ is clear (vectors are added). We identify $\mathfrak{so}(V)$ with $\bigwedge^2 V$; $c \in \bigwedge^2 V$ acts on $\bigwedge V$ by multiplication by $e^c$, and $SO(V)$ in the usual way. Composition in $Gr(A)$ is given by $(v_1, c_1)(v_2, c_2) = (v_1 + v_2, c_1 c_2)$, for $v_{1,2} \in T_x P$, $c_{1,2} \in \bigwedge V$.

We let $\mathbb{R}^*$ act on $\mathcal{T}Cl(A)$: $\lambda \in \mathbb{R}^*$ acts by $\lambda^2$ on $\tau M$; for $\epsilon \neq 0$ it leaves $Cl(A)$ intact and for $\epsilon = 0$ it multiplies $V$ (or $A$) by $\lambda$. Again, there is a natural smooth structure on $\mathcal{T}Cl(A)$ for which this action is smooth.

Finally, the convolution algebra of $Gr(A)$ is canonically isomorphic to the algebra of functions on the supermanifold $\mathcal{E}$: choosing a trivialization of $V$ (a 1st order trivialization, i.e. a connection, is sufficient), this isomorphism is simply the Fourier transform along the fibres of $TM$; obviously (and miraculously as well), it is independent of the choice.

The symbol calculus is ready now. Given an element of the convolution algebra of $\mathcal{T}Cl(A)$, its restriction to $Gr(A)$ is its symbol. A $Cl(A)$-pseudodiff. operator of order $k$ is uniquelly extended to a weight-$k$ $\mathbb{R}^*$-equivariant $\mathcal{T}Cl(A)$-operator.

3

**Letter 7**

The purpose of this note is to put the connection between gerbes and exact CA's into superlanguage. Non-negatively graded supermanifolds (NNGS) will be appearing all the time. They are just supermanifolds with action of the multiplicative semigroup $(\mathbb{R}, \times)$ such that $-1$ acts as the parity operator. Ordinary manifolds are included as NNGS's with trivial action. $Q$-manifolds (in this note) are NNGS's with a degree-one vector field $Q$ with $Q^2 = 0$. (The world of $Q$-manifolds is an amazing place where things like Sullivan's minimal models live.)

The basic example of a $Q$-manifold is $T[1]M$ (here $M$ is a NNGS; notation $T[1]$ is used rather than $\Pi T$ to indicate the grading of the result). If $M$ is a $Q$-manifold, $Q_M$ (as a section of $T[1]M$) gives us a $Q$-equivariant map $M \to T[1]M$. If you like it (guess you don't), $T[1]$ : NNGS's $\to Q$-manifolds is adjoint to the forgetful functor $Q$-manifolds $\to$ NNGS's. If $M$ is some kind of algebroid over $M_0 = 0 \cdot M$, the composition $M \to T[1]M \to T[1]M_0$ is the anchor.

A $Q$-manifold $M$ is a CA if it carries a degree-two symplectic form $\omega$ and a degree-three function $\theta$ with $Q_M = X_\theta$ (a proof of this fairly simple fact is bellow, but I think it should be the definition of CA's). An exact CA is a CA for which the anchor $M \to T[1]M_0$ is a Lagrangian fibration (this is a definition as well).

Now we're getting to gerbes (via a kind of categorification of the vector fields $Q$). Let $\mathbb{R}[n]$ be the 1dim additive group in NNGS's (either $\mathbb{R}$ or $\mathbb{R}^{0|1}$, according to parity of $n$) with the $(\mathbb{R}, \times)$ action $x \mapsto \lambda^n x$. Suppose $N \to M$ is a principal $\mathbb{R}[n]$-bundle (in NNGS's) and moreover $M$ is a $Q$-manifold. We'll construct a principal $\mathbb{R}[n+1]$-bundle $DN \to M$ in $Q$-manifolds. It is done in two steps. The first is a family version of $T[1]$: we have the $Q$-bundle $T[1]N \to T[1]M$, and we also have the $Q$-map $Q_M : M \to T[1]M$; we put $\tilde{N} = Q_M^* T[1]N$. $\tilde{N} \to M$ is already a $Q$-bundle (with structure group $T[1]\mathbb{R}[n]$). Finally (the second step), we put $DN = \tilde{N}/\mathbb{R}[n]$; as promised, it is a principal $(T[1]\mathbb{R}[n])/\mathbb{R}[n] = \mathbb{R}[n+1]$-bundle.

Principal $\mathbb{R}[n]$ bundles in NNGS's are of course trivializable; those in $Q$-manifolds are classified by $n + 1$st cohomology of $Q_M$. These bundles can obviously be added (in both categories; the functor $D$ is clearly additive); this goes to addition in cohomology in the $Q$-case (I wonder if there is some product as well). The $Q$-bundle $DN$ constructed above is trivializable (just because $N$ is) but not canonically. A choice of its trivialization is the same as a choice of a $Q$-structure on $N$ (this is what $DN$ is good for).

Gerbes now appear in this way: Suppose we are given a family $N_i \to U_i \subset M$ of locally defined $N$'s and that the differences $N_i - N_j$ are $Q$-bundles in a coherent way (this is the only instant of a gerbe that we shall need). This means that $DN_i$'s are identified on their overlaps so that we have a well defined global $Q$-bundle. These kinds of gerbes are thus equivalent to locally trivializable principal $\mathbb{R}[n+1]$-bundles in the $Q$-category.

Now we just specialize this abstract nonsense to the simple case when $M = T[1]M_0$ for some ordinary manifold $M_0$. Then everything is locally trivializable. We actually need only $n = 1$. If you wish, a principal $\mathbb{R}[1]$-bundle over $T[1]M_0$ is the same as an extension of vector bundles $0 \to \mathbb{R} \to V \to TM \to 0$; if it's

1

$Q$-equivariant then $V$ is a Lie algebroid. So suppose we are given a family $N_i$ as above, i.e. an infinitesimal DD gerbe. We get as a result a $Q$-equivariant $\mathbb{R}[2]$-principal bundle $K \to T[1]M_0$ (where $K$ is locally $DN_i$). A thing like $K$ is equivalent to an exact CA over $M_0$: just take symplectic reduction of $T^*[2]K$ by $\mathbb{R}[2]$ at moment 1 (this procedure is easily seen to be reversible, but I haven't found any nice geometric way back yet).

Changing the topic completely, we may try to classify symplectic NNGS's with symplectic form of degree $n$. For $n = 1$, $M = T^*[1]M_0$, for $n = 2$ we get Alan's symplectic realization, etc. (cf. below). Then we can take functions $\theta$ of degree $n + 1$ (so that $X_\theta$ is of degree 1) and require $\{\theta, \theta\} = 0$. We get an interesting series starting with Poisson manifolds ($n = 1$) and CA's ($n = 2$). One can also take general $Q$-manifolds to be generalizations of Lie algebroids.

If $(M, \omega)$ is a symplectic NNGS with $\deg \omega = n$ ($n \geq 1$) then $M \to M_0$ is a coisotropic fibration and we can form the corresponding symplectic fibration $M' \to M_0$ by modding out the null directions. The way back form $M'$ to $M$ is easy and unique (this kind of symplectic realization works for bundles (or foliations) of exact symplectic manifolds and is perhaps best visualized via the obvious "contact realization" of a bundle (or foliation) of contact manifolds): the fibres of $M'$ have natural 1-forms $\alpha$ (plug the Louiville vector field into $\omega$) and $M \subset T^*[n]M'$ are simply the covectors which, when restricted to the fibres, are equal to $\alpha$. Now the question is what $M'$ and its fibres may look like. Just take the tangent vector space of a fibre at its origin (the point in $M_0$), notice that it is graded from 1 to $n-1$ and $\omega$ gives you some pairings on the graded components. For $n = 1$ we have $M' = 0$, for $n = 2$ it is $A[1]$ with $A$ a vector bundle with an inner product; for higher $n$'s we can still identify $M'$ with the bundle of tangent spaces at origins (the corresponding groups are contractible), but this is no longer natural (and I don't know what happens for complex manifolds).

# Integration of Courant algebroids and its relatives (a first glance)

**Letter 8**

There is an infinite series of notions that will be denoted as $\Sigma_n$-manifolds and their D-structures (D may mean both Dirac and Dirichlet): symplectic manifolds and their Lagrangian submanifolds ($n = 0$), Poisson manifolds and their coisotropic submanifolds ($n = 1$), Courant algebroids and their (generalized) Dirac structures ($n = 2$); higher $n$'s appear nameless. Poisson manifolds can be (at least locally) integrated to symplectic groupoids; we'll have a look what happens for higher $n$'s. One can say that the result is a symplectic $n$-groupoid. But perhaps more illuminating way of saying almost the same, there will be a symplectic version of $n + 1$ dimensional TFT (i.e. instead of vector spaces we'll have symplectic manifolds, and instead of their elements (or linear maps), Lagrangian submanifolds), accompanied by other structures on boundaries and corners (coming e.g. from D-structures). Though it may not seem obvious at first sight, it is actually very simple. The landscape is vast, however; only few pictures are described here.

As a motivation, it's good to start with integration of Lie algebroids, or even of Lie algebras. Let $A \to M$ be a Lie algebroid; we'd like to construct some groupoid $\Gamma$. Suppose $TI \to A$ is a Lie algebroid morphism (covering a map $I \to M$, where $I$ is an interval); it certainly gives us a morphism in $\Gamma$ between the endpoints $P$ and $Q$ (really obtained by composing infinitesimal morphisms along $I$):

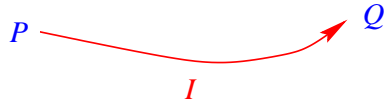

However, there are many ways how to obtain a given morphism in $\Gamma$ using this construction. Consider now this picture:

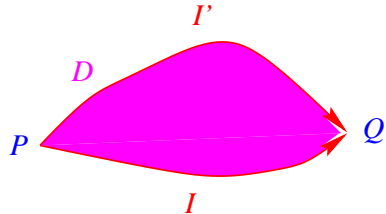

$D$ is a disk and we consider a Lie algebroid morphism $TD \to A$. In this case the morphisms in $\Gamma$ given by $I$ and $I'$ are certainly equal. Locally, these are the only identifications we have to make (of course, global smooth groupoid needn't exist and if it does, it is not specified completely by $A$ either; these questions are *not* addressed in this elementary note); we've integrated $A$ to a (at least locally smooth) groupoid. It looks very much like Poincaré groupoid.

Now some terminology: an *N-manifold* (shorthand for non-negatively graded supermanifold) is a supermanifold with action of the multiplicative semigroup $(\mathbb{R}, \times)$ such that $-1$ acts as the parity operator. An *NQ-manifold* is an N-manifold with a degree-one vector field $Q$ of square 0. A $\Sigma_n$-*manifold* is an NQ-manifold with $Q$-invariant symplectic form of degree $n$; finally, a *D-structure* is its Lagrangian NQ-submanifold.

A basic example of an NQ-manifold is $T[1]M$, where $M$ is an ordinary manifold (the notation $T[1]$ is used rather than $\Pi T$ to indicate the grading). NQ-manifolds should be understood as generalizations of this, i.e. as 'differential models' of manifolds. Likewise, NQ-maps from $T[1]M$ to an NQ-manifold are analogous to maps from $M$ to an ordinary manifold; this principle is used when we speak about homotopy of NQ-manifolds.

Lie algebroids are special cases of NQ-manifolds. Above, we considered NQ-maps from $T[1]I$ and $T[1]D$. Thus, the global object (groupoid) was connected with homotopy properties of the NQ-manifold - it was really the Poincaré groupoid. This will be the general principle: we shall consider maps from $T[1]M$ ($M$ is a manifold or perhaps a simplicial complex) to our NQ-manifold and shall watch what will happen.

NQ-manifolds can be seen as generalizations of Lie algebroids, as they can be integrated to $n$-groupoids. First we'll describe a measure of complexity of an NQ-manifold that may be called its *degree*. Namely, let $X$ be an NQ-manifold and $P$ a point in $0 \cdot X$ (this 0 is from the semigroup $(\mathbb{R}, \times)$). The semigroup acts on $T_P X$; the highest weight that appears is the degree at $P$. It is locally constant in $0 \cdot X$. In particular, degree 1 means that $X$ is a Lie algebroid.

Now an important thing to notice about NQ-manifolds is that they can have nontrivial local homotopy up to the degree (i.e. small maps from $T[1]S^k$ can in general be extended to small maps from $T[1]B^{k+1}$ only if $k$ is larger than the degree). We have to use *Poincaré n-groupoids*: objects are points, 1-morphisms paths between points, 2-morphisms disks between two paths connecting the same two points, etc., only $n$-morphisms are homotopy classes of $n$-balls. As we see, a good choice is $n = $ degree. Thus, $X$ is integrated to an $n$-groupoid with objects in $0 \cdot X$. Vice versa, $X$ can be recovered from this $n$-groupoid.

As an example, here are NQ-manifolds connected with $U(1)$-bundles, $U(1)$-gerbes, etc., or more precisely, with their gauge groupoids (as you expect, the 'gauge groupoid' of a gerbe is actually a 2-groupoid, etc.). Let $\mathbb{R}[n]$ denote the additive group in degree $n$ (i.e. $\mathbb{R}[n]$ is either $\mathbb{R}$ or $\mathbb{R}^{0|1}$ as a supergroup, according to the parity of $n$, with the $(\mathbb{R}, \times)$ action $x \mapsto \lambda^n x$, and with zero $Q$-field). Let $M$ be an ordinary manifold; we shall consider principal $\mathbb{R}[n]$-bundles $X \to T[1]M$ (both the projection $X \to T[1]M$ and the action $\mathbb{R}[n] \times X \to X$ are NQ-maps). Such $X$'s are classified by $H^{n+1}(M, \mathbb{R})$ and their degree is $n$ (unless $n = 0$, when it is 1). Let us consider the $n$-groupoid corresponding to such an $X$. As it turns out, the objects are points in $M$, 1-morphisms paths in $M$ between the points, etc.; however, $n$-morphisms are not just homotopy classes of $n$-balls in $M$, but instead they are from a central extension of them by $\mathbb{R}$. Ultimately, this $\mathbb{R}$ is to be changed to $U(1)$ and $H^{n+1}(M, \mathbb{R})$ to $H^{n+1}(M, \mathbb{Z})$, but as I said before, these questions are not considered here.

This is all I can say about general NQ-manifolds. Now we pass to $\Sigma_n$-manifolds. The degree of a $\Sigma_n$ manifold is $n$ (it can be less if $0 \cdot X$ is a single point). There is a way how to say that the corresponding $n$-groupoid is symplectic, and vice versa, the symplectic structure on the $n$-groupoid makes the NQ-manifold into a $\Sigma_n$-manifold. The definition is not attempted here (but some of the flavor is in what follows). Instead we shall consider symplectic TFT's with boundaries coloured by D-structures.

Let $M$ be a closed oriented $n$-dimensional manifold and let $X$ be a $\Sigma_n$-manifold. The space of maps $T[1]M \to X$ carries a closed 2-form; this form is symplectic on the space of all maps, but we shall consider only the space of NQ-maps. As it turns out, modding out the null directions of the 2-form is equivalent to passing to homotopy classes of maps. The 2-form is defined as follows: given two infinitesimal deformations $u$ and $v$ of a map $T[1]M \to X$, the value of the form on them is $\int_{T[1]M} \omega_X(u, v)$. For example, if $n = 2$ and $X = \mathfrak{g}[1]$ where $\mathfrak{g}$ is a Lie algebra with invariant inner product, we get the moduli space of flat $\mathfrak{g}$-connections (the symplectic 2-groupoid is also well known in this case: 2-morphisms are loops in the simply-connected $G$ modulo right (say) translations (since $\pi_2(G) = 0$); the symplectic form on their space is the famous one).

It is also useful to choose several D-structures in $X$ and to consider $M$'s with

boundaries coloured in variuos colours, each colour corresponding to one of the D-structures. For example, we can get the double symplectic groupoid corresponding to a Lie bialgebra or even to a Lie bialgebroid:

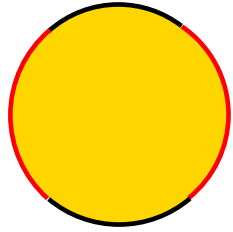

The condition on the boundary of the disk is that each colour (red or black) on the boundary of the disk is mapped to one of the two transversal D-structures. It is interesting to notice how close the double groupoid is to the moduli space of flat connections.

We got symplectic spaces for $n$-dim manifolds, perhaps with coloured boundaries; moreover, the NQ maps from $T[1]$ of an $n+1$-dim manifold give us a Lagrangian submanifold in the symplectic space of its boundary, i.e. we have a symplectic version of TFT. For example,

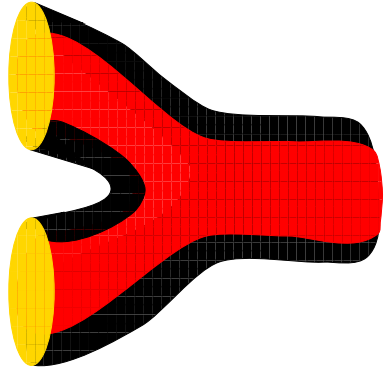

gives one of the two products in the double symplectic groupoid. The remaining operations (together with reasons why they are compatible in the appropriate way) are in hep-th/9811136.

There is one more reason to be interested in $\Sigma_n$-manifolds. They are connected with variational problems: extremal $n$-folds are just $n$-dim submanifolds of a $\Sigma_n$-manifold $X$ (i.e. NQ maps $T[1]M \to X$ with $n$-dim $M$); $X$ encodes the type of the system, while Lagrangian appears as an additional constraint on the map (D-structures play the role of boundary contitions). This may lead to very nice applications; up to now, I have only partially explored the case $n = 2$ (Poisson-Lie T-duality). When $\Sigma_n$-manifolds are eventually quantized (via $n+1$-dim TFT's), this may lead to some quantum applications as well (QFT's on boundary of TFT's).

*References:* This note is an almost immediate consequence of two papers. One is Sullivan's Infinitesimal computations. I used only few ideas – the philosophy of NQ-manifolds as differential models and the idea of integration using homotopy properties. The paper contains much more, so more can be expected. The second main reference is AKSZ's paper on BV quantization and TFT's. All the symplectic ideas of this note come from there. The paper gives direct instructions how to quantize $\Sigma_n$-manifolds using perturbative TFT (instead of the rather naive symplectic TFT described above). Poincaré $n$-groupoids are (I think) from the 600pp Grothendieck's postcard to Quillen (unfortunately, I haven't seen it yet). The recent construction of symplectic groupoids by Cattaneo and Felder was also very useful.


 
\end{document}